\documentclass[11pt]{article}
\usepackage{amssymb}

\usepackage{amsmath}
\usepackage{theorem}
\usepackage{here}
\usepackage[dvipdfmx]{color}

\newtheorem{thm}{Theorem}[section]

\newtheorem{prp}{Proposition}[section]
\theorembodyfont{\rmfamily}

\makeatletter

\@addtoreset{equation}{section}
\makeatother

\usepackage{comment} %
\usepackage{bm}
\usepackage[dvipdfmx]{graphicx} %
\def\Kt{{\widetilde K}}

\def\Th{{\Theta}}
\def\De{{\Delta}}

\def\d{{\text{\boldmath $d$}}}

\def\x{{\text{\boldmath $x$}}}

\def\K{{\text{\boldmath $K$}}}

\def\X{{\text{\boldmath $X$}}}
\def\Y{{\text{\boldmath $Y$}}}

\def\Kh{{\widehat K}}

\def\Kbh{{\widehat \K}}
\def\Kbt{{\widetilde \K}}
\def\[{{\text{\boldmath $[$}}}
\def\]{{\text{\boldmath $]$}}}

\def\/{{\Bigr/\!\!}}

\def\1r{{\rm (1)}}
\def\2r{{\rm (2)}}
\def\3r{{\rm (3)}}
\def\4r{{\rm (4)}}
\def\5r{{\rm (5)}}

\def\non{{\nonumber}}
\begin{document}
\title{Bayesian Estimation for the Multivariate Hypergeometric Distribution Incorporating Information from Aggregated Observations\footnote{This preprint has not undergone peer review (when applicable) or any post-submission improvements or corrections. 
The Version of Record of this article is published in Japanese Journal of Statistics and Data Science, and is available online at https://doi.org/10.1007/s42081-023-00224-z. }}
\author{
Yasuyuki Hamura\footnote{Graduate School of Economics, Kyoto University, 
Yoshida-Honmachi, Sakyo-ku, Kyoto, 606-8501, JAPAN. 
\newline{
E-Mail: yasu.stat@gmail.com}} \
}
\maketitle
\begin{abstract}
In this short note, we consider the problem of estimating multivariate hypergeometric parameters under squared error loss when side information in aggregated data is available. 
We use the symmetric multinomial prior to obtain Bayes estimators. 
It is shown that by incorporating the side information, we can construct an improved estimator. 

\par\vspace{4mm}
\end{abstract}

\section{Introduction}
\label{sec:introduction}
The multivariate hypergeometric distribution is a multivariate generalization of the usual (univariate) hypergeometric distribution and arises from sampling without replacement when there are multiple categories (see, for example, Johnson, Kotz and Balakrishnan (1997)). 
The multivariate hypergeometric distribution will thus be one of the most basic discrete distributions. 
For some applications of the distribution,  see, for example, %
Johnson, Kotz and Wu (1991). 
Although some authors (e.g., Hodges and Lehmann (1950), Trybula (1958), and Wilczynski (1985)) considered this distribution in the context of decision theory, relatively little attention has been paid to it compared with more famous discrete distributions such as the Poisson and multinomial distributions. 

In this note, we construct an improved estimator for multivariate hypergeometric parameters in a framework similar to those of Hamura (2021, 2022). 
Specifically, we assume that a finite population of $K \in \mathbb{N}$ individuals is composed of $m \in \mathbb{N}$ mutually exclusive subgroups and that each subgroup is divided into $n \in \mathbb{N}$ categories. 
We consider the problem of simultaneously estimating the numbers of the elements of all these $m n $ categories, $(( K_{i, j} )_{j = 1}^{n} )_{i = 1}^{m}$, $\sum_{i = 1}^{m} \sum_{j = 1}^{n} K_{i, j} = K$, under squared error loss. 
We suppose that we first sample $L \in \mathbb{N}$ individuals without replacement to observe a multivariate hypergeometric vector $\X = (( X_{i, j} )_{j = 1}^{n} )_{i = 1}^{m}$, where $X_{i, j}$ is the number of individuals sampled from the $(i, j)$-th category, $j = 1, \dots , n$, $i = 1, \dots , m$, and that we then sample, without replacement, $L' \in \mathbb{N}$ individuals from the remaining population of size $K - L$ but for some reason observe $\Y = \big( \sum_{j = 1}^{n} Y_{i, j}^{*} \big) _{i = 1}^{m}$ only, where for each $i = 1, \dots , m$ and $j = 1, \dots , n$, $Y_{i, j}^{*}$ denotes the number of individuals sampled in this second survey from the $(i, j)$-th category having $K_{i, j} - X_{i, j}$ remaining individuals. 
Following Balakrishnan and Ma (1996), we assume a multinomial prior for $\K = (( K_{i, j} )_{j = 1}^{n} )_{i = 1}^{m}$ to derive Bayesian estimators. 
We investigate when using both $\X $ and $\Y $ is better than using $\X $ only.

\section{The Problem}
\label{sec:problem} 
Let $K \in \mathbb{N}$ be a given number. 
Let $m, n \in \mathbb{N}$ be known numbers and let $(( K_{i, j} )_{j = 1}^{n} )_{i = 1}^{m} \in \Th = \big\{ (( \mathring{K} _{i, j} )_{j = 1}^{n} )_{i = 1}^{m} \in {\mathbb{N} _0}^{m n} \big| \sum_{i = 1}^{m} \sum_{j = 1}^{n} \mathring{K} _{i, j} = K \big\} $ be unknown parameters, where $\mathbb{N} _0 = \{ 0 \} \cup \mathbb{N} = \{ 0, 1, 2, \dotsc \} $. 
Suppose that $\X = (( X_{i, j} )_{j = 1}^{n} )_{i = 1}^{m} \in {\mathbb{N} _0}^{m n}$ is distributed with mass function 
\begin{align}
\Big( \prod_{i = 1}^{m} \prod_{j = 1}^{n} \binom{K_{i, j}}{x_{i, j}} \Big) / \binom{K}{L} \text{,} \non 
\end{align}
$(( x_{i, j} )_{j = 1}^{n} )_{i = 1}^{m} \in D = \big\{ (( \mathring{x} _{i, j} )_{j = 1}^{n} )_{i = 1}^{m} \in {\mathbb{N} _0}^{m n} \big| \sum_{i = 1}^{m} \sum_{j = 1}^{n} \mathring{x} _{i, j} = L \big\} \cap \prod_{i = 1}^{m} \prod_{j = 1}^{n} [0, K_{i, j} ]$, for known $L \in \mathbb{N} \cap [0, K]$. 
Suppose further that for known $L' \in \mathbb{N} \cap [0, K - L]$, the random vector $\Y = ( Y_i )_{i = 1}^{m} \in {\mathbb{N} _0}^m$ is such that for each $\x = (( x_{i, j} )_{j = 1}^{n} )_{i = 1}^{m} \in D$, the conditional mass function of $\Y | ( \X = \x )$ is given by 
\begin{align}
\Big( \prod_{i = 1}^{m} \binom{K_{i, \cdot } - x_{i, \cdot }}{y_i} \Big) / \binom{K - L}{L'} \text{,} \non 
\end{align}
$( y_i )_{i = 1}^{m} \in \big\{ ( \mathring{y} )_{i = 1}^{m} \in {\mathbb{N} _0}^m \big| \sum_{i = 1}^{m} \mathring{y} _i = L' \big\} \cap \prod_{i = 1}^{m} [0, K_{i, \cdot } - x_{i, \cdot } ]$, where $K_{i, \cdot } = \sum_{j = 1}^{n} K_{i, j}$ and $x_{i, \cdot } = \sum_{j = 1}^{n} x_{i, j}$ for $i = 1, \dots , m$. 
We consider the estimation of $\K = (( K_{i, j} )_{j = 1}^{n} )_{i = 1}^{m}$ on the basis of $\X $ and $\Y $ under the squared error loss 
\begin{align}
L( \d , \K ) &= \sum_{i = 1}^{m} \sum_{j = 1}^{n} ( d_{i, j} - K_{i, j} )^2 \text{,} \non 
\end{align}
$\d = (( d_{i, j} )_{j = 1}^{n} )_{i = 1}^{m} \in \mathbb{R} ^{m n}$.

\section{Bayesian Estimators}
\label{sec:estimator} 
We use the multinomial prior 
\begin{align}
\K \sim \pi ( \K ) = {K ! \over \prod_{i = 1}^{m} \prod_{j = 1}^{n} ( K_{i, j} !)} \prod_{i = 1}^{m} \prod_{j = 1}^{n} {p_{i, j}}^{K_{i, j}} \text{,} \non 
\end{align}
where $(( p_{i, j} )_{j = 1}^{n} )_{i = 1}^{m} \in \big\{ (( \mathring{p} _{i, j} )_{j = 1}^{n} )_{i = 1}^{m} \in (0, 1)^{m n} \big| \sum_{i = 1}^{m} \sum_{j = 1}^{n} \mathring{p} _{i, j} = 1 \big\} $. 
Let $p_{i, \cdot } = \sum_{j = 1}^{n} p_{i, j}$ for $i = 1, \dots , m$. 

\begin{prp}
\label{prp:estimator} 
\hfill
\begin{itemize}
\item[{\rm{(i)}}]
The Bayes estimator of $\K $ with respect to the above prior and the observations $\X $ and $\Y $ is 
\begin{align}
\Kbh ( \X , \Y ) = (( \Kh _{i, j} ( \X , \Y ))_{j = 1}^{n} )_{i = 1}^{m} = \Big( \Big( X_{i, j} + {p_{i, j} \over p_{i, \cdot }} Y_i + p_{i, j} (K - L - L' ) \Big) _{j = 1}^{n} \Big) _{i = 1}^{m} \text{.} \non 
\end{align}
\item[{\rm{(ii)}}]
The Bayes estimator of $\K $ with respect to the above prior and $\X $ only is 
\begin{align}
\Kbh ( \X )= (( \Kh _{i, j} ( \X ))_{j = 1}^{n} )_{i = 1}^{m} = (( X_{i, j} + p_{i, j} (K - L))_{j = 1}^{n} )_{i = 1}^{m} \text{.} \non 
\end{align}
\end{itemize}
\end{prp}

\noindent
{\bf Proof%
.} \ \ Let $X_{i, \cdot } = \sum_{j = 1}^{n} X_{i, j}$ for $i = 1, \dots , m$. 
Then the posterior mass of $\K $ is 
\begin{align}
p( \K | \X , \Y ) &\propto \Big( {1 \over \prod_{i = 1}^{m} \prod_{j = 1}^{n} ( K_{i, j} !)} \prod_{i = 1}^{m} \prod_{j = 1}^{n} {p_{i, j}}^{K_{i, j}} \Big) \Big( \prod_{i = 1}^{m} \binom{K_{i, \cdot } - X_{i, \cdot }}{Y_i} \Big) \prod_{i = 1}^{m} \prod_{j = 1}^{n} \binom{K_{i, j}}{X_{i, j}} \non \\
&\propto \prod_{i = 1}^{m} \Big\{ {1 \over ( K_{i, \cdot } - X_{i, \cdot } - Y_i ) !} {p_{i, \cdot }}^{K_{i, \cdot }} {( K_{i, \cdot } - X_{i, \cdot } ) ! \over \prod_{j = 1}^{n} ( K_{i, j} - X_{i, j} ) !} \prod_{j = 1}^{n} \Big( {p_{i, j} \over p_{i, \cdot }} \Big) ^{K_{i, j}} \Big\} \non \\
&\propto \Big[ \prod_{i = 1}^{m} \Big\{ {1 \over ( K_{i, \cdot } - X_{i, \cdot } - Y_i ) !} {p_{i, \cdot }}^{K_{i, \cdot } - X_{i, \cdot } - Y_i} \Big\} \Big] \prod_{i = 1}^{m} \Big\{ {( K_{i, \cdot } - X_{i, \cdot } ) ! \over \prod_{j = 1}^{n} ( K_{i, j} - X_{i, j} ) !} \prod_{j = 1}^{n} \Big( {p_{i, j} \over p_{i, \cdot }} \Big) ^{K_{i, j} - X_{i, j}} \Big\} \text{,} \non 
\end{align}
where $1 / (t !) = 0$ for any $t < 0$. 
Let $\Kt _i = K_{i, \cdot } - X_{i, \cdot } - Y_i$ for $i = 1, \dots , m$ and let $\Kt _{i, j} = K_{i, j} - X_{i, j}$ for $j = 1, \dots , n$ for $i = 1, \dots m$. 
Then the posterior mass of $\Kbt = ((( \Kt _{i, j} )_{j = 1}^{n} )_{i = 1}^{m} , ( \Kt _i )_{i = 1}^{m} )$ is 
\begin{align}
p( \Kbt | \X , \Y ) &\propto {\rm{Multin}}_{m - 1} (( \Kt _i )_{i = 1}^{m} | K - L - L' , ( p_{i, \cdot } )_{i = 1}^{m} ) \prod_{i = 1}^{m} {\rm{Multin}}_{n - 1} (( \Kt _{i, j} )_{j = 1}^{n} | Y_i + \Kt _i , ( p_{i, j} / p_{i, \cdot } )_{j = 1}^{n} ) \text{.} \non 
\end{align}
Therefore, 
\begin{align}
E[ K_{i, j} | \X , \Y ] &= E[ X_{i, j} + \Kt _{i, j} | \X , \Y ] \non \\
&= X_{i, j} + E[ ( Y_i + \Kt _i ) p_{i, j} / p_{i, \cdot } | \X , \Y ] \non \\
&= X_{i, j} + \{ Y_i + (K - L - L' ) p_{i, \cdot } \} p_{i, j} / p_{i, \cdot } \non 
\end{align}
for all $j = 1, \dots , n$ for all $i = 1, \dots , m$, which is the desired result. 
\hfill$\Box$

\section{Dominance}
\label{subsec:dominance} 
Here, we put $p_{i, j} = 1 / (m n)$ for all $j = 1, \dots , n$ for all $i = 1, \dots , m$ and compare the risk functions of 
\begin{align}
\Kbh ( \X , \Y ) = (( \Kh _{i, j} ( \X , \Y ))_{j = 1}^{n} )_{i = 1}^{m} = \Big( \Big( X_{i, j} + {Y_i \over n} + {K - L - L' \over m n} \Big) _{j = 1}^{n} \Big) _{i = 1}^{m} \text{.} \non 
\end{align}
and 
\begin{align}
\Kbh ( \X ) = (( \Kh _{i, j} ( \X ))_{j = 1}^{n} )_{i = 1}^{m} = \Big( \Big( X_{i, j} + {K - L \over m n} \Big) _{j = 1}^{n} \Big) _{i = 1}^{m} \text{.} \non 
\end{align}

\begin{thm}
\label{thm:dominance} 
Suppose that $m \ge 2$. 
Suppose that $L \le K - 2$ and that $L' \ge 1$. 
Then $\Kbh ( \X , \Y )$ dominates $\Kbh ( \X )$ if $2 L + L' \ge K$. 
Furthermore, if $K / m \in \mathbb{N}$, then $\Kbh ( \X , \Y )$ dominates $\Kbh ( \X )$ if and only if $2 L + L' \ge K$. 
\end{thm}

Although the condition $2 L + L' \ge K$ is restrictive, it is necessary for $\Kbh ( \X , \Y )$ to dominate $\Kbh ( \X )$ when $K / m \in \mathbb{N}$. 
It follows from the following proof that if $K_{1, \cdot } = \dots = K_{m, \cdot }$, then we have $E[ L( \Kbh ( \X , \Y ), \K ) ] - E[ L( \Kbh ( \X ), \K ) ] = \max_{\mathring{\K } \in \Th } \{ E[ L( \Kbh ( \X , \Y ), \mathring{\K } ) ] - E[ L( \Kbh ( \X ), \mathring{\K } ) ] \} $. 

\bigskip

\noindent
{\bf Proof of Theorem \ref{thm:dominance}.} \ \ Fix $i = 1, \dots , m$ and $j = 1, \dots , n$. 
Then 
\begin{align}
E[ \{ \Kh _{i, j} ( \X , \Y ) - K_{i, j} \} ^2 ] &= E[ V( Y_i | \X ) / n^2 + \{ E[ Y_i | \X ] / n + X_{i, j} + (K - L - L' ) / (m n) - K_{i, j} \} ^2 ] \non \\
&= E \Big[ {L' (K - L - L' ) \over (K - L)^2 (K - L - 1)} ( K_{i, \cdot } - X_{i, \cdot } ) \{ K - L - ( K_{i, \cdot } - X_{i, \cdot } ) \} / n^2 \non \\
&\quad + \Big\{ {L' \over K - L} ( K_{i, \cdot } - X_{i, \cdot } ) / n + X_{i, j} + (K - L - L' ) / (m n) - K_{i, j} \Big\} ^2 \Big] \text{,} \non 
\end{align}
where the second equality follows from (39.7) and (39.8) of Johnson, Kotz and Balakrishnan (1997). 
Meanwhile, 
\begin{align}
E[ \{ \Kh _{i, j} ( \X ) - K_{i, j} \} ^2 ] &= %
E[ \{ X_{i, j} + (K - L - L' ) / (m n) - K_{i, j} + L' / (m n) \} ^2 ] \text{.} \non 
\end{align}
Therefore, 
\begin{align}
&E[ \{ \Kh _{i, j} ( \X , \Y ) - K_{i, j} \} ^2 ] - E[ \{ \Kh _{i, j} ( \X ) - K_{i, j} \} ^2 ] \non \\
&= E \Big[ {L' (K - L - L' ) \over (K - L)^2 (K - L - 1)} ( K_{i, \cdot } - X_{i, \cdot } ) \{ K - L - ( K_{i, \cdot } - X_{i, \cdot } ) \} / n^2 \non \\
&\quad + \Big\{ {L' \over K - L} ( K_{i, \cdot } - X_{i, \cdot } ) / n - L' / (m n) \Big\} \non \\
&\quad \times \Big[ 2 \{ X_{i, j} + (K - L - L' ) / (m n) - K_{i, j} \} + {L' \over K - L} ( K_{i, \cdot } - X_{i, \cdot } ) / n + L' / (m n) \Big] \Big] \text{.} \label{tdominancep1} 
\end{align}

Let $\De = E[ L( \Kbh ( \X , \Y ), \K ) ] - E[ L( \Kbh ( \X ), \K ) ]$. 
Then, by (\ref{tdominancep1}), 
\begin{align}
\De &= \sum_{i = 1}^{m} E \Big[ {L' (K - L - L' ) \over (K - L)^2 (K - L - 1)} ( K_{i, \cdot } - X_{i, \cdot } ) \{ K - L - ( K_i - X_{i, \cdot } ) \} / n \non \\
&\quad + \Big\{ {L' \over K - L} ( K_{i, \cdot } - X_{i, \cdot } ) / n - L' / (m n) \Big\} \non \\
&\quad \times \Big[ 2 \{ X_{i, \cdot } + (K - L - L' ) / m - K_{i, \cdot } \} + {L' \over K - L} ( K_{i, \cdot } - X_{i, \cdot } ) + L' / m \Big] \Big] \non \\
&= {L' \over n} \sum_{i = 1}^{m} E \Big[ {K - L - L' \over K - L - 1} {K_{i, \cdot } - X_{i, \cdot } \over K - L} \Big( 1 - {K_i - X_{i, \cdot } \over K - L} \Big) \non \\
&\quad + \Big( {K_{i, \cdot } - X_{i, \cdot } \over K - L} - 1 / m \Big) \non \\
&\quad \times \Big\{ 2 (K - L) \Big( - {K_{i, \cdot } - X_{i, \cdot } \over K - L} + {K - L - L' \over K - L} / m \Big) + {K_{i, \cdot } - X_{i, \cdot } \over K - L} L' + L' / m \Big\} \Big] \non \\
&= {L' \over n} \sum_{i = 1}^{m} E \Big[ {K - L - L' \over K - L - 1} {K_{i, \cdot } - X_{i, \cdot } \over K - L} - {K - L - L' \over K - L - 1} \Big( {K_{i, \cdot } - X_{i, \cdot } \over K - L} \Big) ^2 \non \\
&\quad + \Big( {K_{i, \cdot } - X_{i, \cdot } \over K - L} - 1 / m \Big) \non \\
&\quad \times \Big[ - \{ 2 (K - L) - L' \} {K_{i, \cdot } - X_{i, \cdot } \over K - L} + 2 (K - L) {K - L - L' \over K - L} / m + L' / m \Big] \Big] \non 
\end{align}
Therefore, 
\begin{align}
\De &= {L' \over n} E \Big[ {K - L - L' \over K - L - 1} - {K - L - L' \over K - L - 1} \sum_{i = 1}^{m} \Big( {K_{i, \cdot } - X_{i, \cdot } \over K - L} \Big) ^2 \non \\
&\quad - \{ 2 (K - L) - L' \} \sum_{i = 1}^{m} \Big( {K_{i, \cdot } - X_{i, \cdot } \over K - L} \Big) ^2 + 2 (K - L) {K - L - L' \over K - L} / m + L' / m \non \\
&\quad - (1 / m) \Big[ - \{ 2 (K - L) - L' \} + 2 (K - L) {K - L - L' \over K - L} + L' \Big] \Big] \non \\
&= {L' \over n} E \Big[ - \Big\{ 2 (K - L) - L' + {K - L - L' \over K - L - 1} \Big\} \sum_{i = 1}^{m} \Big( {K_{i, \cdot } - X_{i, \cdot } \over K - L} \Big) ^2 + (1 / m) \{ 2 (K - L) - L' \} + {K - L - L' \over K - L - 1} \Big] \text{.} \non 
\end{align}
Since 
\begin{align}
E[ ( K_{i, \cdot } - X_{i, \cdot } )^2 ] &= V( X_{i, \cdot } ) + (E[ X_{i, \cdot } ] - K_{i, \cdot } )^2 \non \\
&= {L (K - L) \over K^2 (K - 1)} K_{i, \cdot } (K - K_{i, \cdot } ) + \Big( {L \over K} K_{i, \cdot } - K_{i, \cdot } \Big) ^2 \non \\
&= {(K - L) \{ (K - 1) (K - L) - L \} \over K^2 (K - 1)} {K_{i, \cdot }}^2 + {L (K - L) \over K (K - 1)} K_{i, \cdot } \non 
\end{align}
for all $i = 1, \dots , m$ by the fact that $( X_{i, \cdot } )_{i = 1}^{m}$ is also multivariate hypergeometric, we have 
\begin{align}
\sum_{i = 1}^{m} E[ ( K_{i, \cdot } - X_{i, \cdot } )^2 ] &\ge m \Big\{ {(K - L) (K - L - 1) \over K (K - 1)} {K^2 \over m^2} + {L (K - L) \over K (K - 1)} {K \over m} \Big\} \non \\
&= {(K - L) (K - L - 1) \over K - 1} {K \over m} + {L (K - L) \over K - 1} \non 
\end{align}
by Jensen's inequality. 
The above inequality is strict if and only if $K_{1, \cdot } , \dots , K_{m, \cdot }$ are not all equal. 
Thus, %
\begin{align}
\De &\le %
{L' \over n} \Big[ - \Big\{ 2 (K - L) - L' + {K - L - L' \over K - L - 1} \Big\} {1 \over (K - L)^2} \Big\{ {(K - L) (K - L - 1) \over K - 1} {K \over m} + {L (K - L) \over K - 1} \Big\} \non \\
&\quad + (1 / m) \{ 2 (K - L) - L' \} + {K - L - L' \over K - L - 1} \Big] \text{,} \non %
\end{align}
the right-hand side of which is equal to 
\begin{align}
{L' \over n} \Big( 1 - {1 \over m} \Big) {K - 2 L - L' \over K - 1} \text{,} \non 
\end{align}
and this completes the proof. 
\hfill$\Box$

\section{Discussion}
\label{sec:discussion} 
In this note, we derived one dominance result for one multivariate hypergeometric model involving aggregated observations. 
There are many other important models involving the multivariate hypergeometric distribution. 
For example, we could first observe aggregated data and then observe unaggregated data. 
Aggregated and unaggregated data could be independently observed. 
Unaggregated multivariate hypergeometric data could be gathered from each subgroup. 
Although Bayesian estimators are more complicated in general, these cases will certainly be worth considering from a practical point of view.

\end{document}